\documentclass[12pt]{article}
\usepackage[english]{babel}
\usepackage{amsmath,amsfonts,amssymb,amsthm} 

\usepackage{mathdots}

\usepackage{epsf}
\usepackage{graphicx} 

\usepackage[numbers]{natbib}

\usepackage{subfigure}
\usepackage{float}
\usepackage{tikz}
\usetikzlibrary{decorations.pathreplacing,positioning,calc,patterns,patterns.meta}

\usepackage{url}
\usepackage{xcolor}
\usepackage{hyperref} 

\usepackage{makecell,diagbox}

\newcommand\restr[2]{{
  \left.\kern-\nulldelimiterspace 
  #1 
  \vphantom{\big|} 
  \right|_{#2} 
  }} 

\title{Towards modular Hierarchical Poincaré-Steklov solvers}
\author{Michal Outrata and Jos{\'{e}} Pablo Lucero Lorca}
\date{\vspace{-5ex}}

\begin{document}
\maketitle

\begin{abstract}
We revisit the Hierarchical Poincaré–Steklov (HPS) method for the Poisson equation using standard Q1 finite elements, building on the original in \cite{Martinsson2013}. While corner degrees of freedom were implicitly handled in that work, subsequent spectral-element implementations have typically avoided them. In Q1-FEM, however, corner coupling cannot be factored out, and we show how the HPS merge procedure naturally accommodates it when corners are enclosed by elements. This clarification bridges a conceptual gap between algebraic Schur-complement methods and operator-based formulations, providing a consistent path for the FEM community to adopt HPS to retain the Poincaré–Steklov interpretation at both continuous and discrete levels.
\end{abstract}

\section{Introduction}

Hierarchical Poincaré--Steklov (HPS) solvers are a class of hierarchical direct solvers designed for elliptic PDEs; the name was coined in~\cite{GillmanMartinsson2014,Martinsson2015} but some of the ideas are in~\cite{Chen2002, MartinssonRokhlin2005,  GreengardGueyffierMartinssonRokhlin2009, Gillman2011}. Starting with many subdomains, the goal is to recursively merge local boundary operators -- typically \emph{Dirichlet-to-Neumann} (DtN) or \emph{Impedance-to-Impedance} (ItI) maps -- constructing a \emph{global} one that we apply to the problem. The HPS approach achieves high accuracy and near-optimal complexity, combining ideas present in hierarchical matrix computations ($\mathcal{H}$), domain decomposition (DD) and direct solvers, i.e., it is poised to be the keystone connecting several communities, e.g., BDDC or FETI; in words of one of the authors after reading~\cite{lu2024compression}: ``For the love of God, they need to start talking to each other!''. In our opinion this communication has been limited also due to the strong spectral element methods (SEM) background of the HPS community; the formulation, discretization and computation in HPS are often entangled together, making it difficult to relate pros and cons of the ``package'' to its parts. We want to provide a \emph{modular} alternative, approachable for readers across multiple communities and for the sake of space we focus on the \emph{corner points}\footnote{We note that in the DD community we usually use the term \emph{cross points}, e.g.,~\cite{gander2014cross,ClaeysParolin2022}.}, which are \emph{routinely} considered an obstacle in the HPS context~\cite{GillmanMartinsson2014,Martinsson2015,GillmanBarnettMartinsson2015, BabbGillmanHaoMartinsson2018,BeamsGillmanHewett2020,LuceroLorcaBeamsBeecroftGillman2024}. Many aspects of what follows can be found \emph{somewhere} in the literature, sometimes with limited references to the other fields but, to the best of our knowledge, a modular HPS exposition is \emph{nowhere to be found} in the HPS literature. 

As mentioned, the HPS community is using predominantly SEM on tensor product grids -- it offers high (possible) accuracy \emph{and} lets us avoid the \emph{corner points}, e.g., with the Gauss-Legendre points. If the corner points appear, in SEM they usually come \emph{decoupled} from the interior nodes or can be avoided altogether by modifying the spectral discretization, see, e.g.,~\cite{GillmanMartinsson2014, GillmanBarnettMartinsson2015,BabbGillmanHaoMartinsson2018,BeamsGillmanHewett2020}. HPS using finite differences or finite volumes,~\cite{Gillman2011,ChipmanCalhounBurstedde2024}, also rely on avoiding the corner-coupling issue that arises in, e.g., FEM. The rationale is both \emph{analytical} but also practical: the used Poincar{\'{e}}-Steklov (PS) operators need not be well-defined in the presence of corners and the tensor-product basis naturally isolates corner DoFs~\cite{GeldermansGillman2019, Fortunato2024, MeliaFortunatoGillmanONeil2025, KumpYesypenkoMartinsson2025}. Hence, for many new readers, the HPS methods are intrinsically connected with such discretization schemes.

However, essentially the same problems have been studied also from the algebraic perspective, e.g., nested dissection, hierarchical semi-separable and hierarchical multifrontal techniques, e.g.,~\cite{george1973nested,GrasedyckKriemannLeBorne2009, XiaChandrasekaranGuLi2009,SchmitzYing2012}, are \emph{purely algebraic}: they operate directly on the discrete system, exploiting observed numerical blockwise low-rankness for compression and factorization. Although the foundational work on hierarchical matrices, see e.g.,~\cite{borm2006hierarchicalmatrices,hackbusch2015hierarchical}, is built on the continuous operators, to the best of our knowledge, it does not include PS operators, nor incorporate static condensation or skeletonization. The recursive skeletonization can be viewed within the multilevel DD or multigrid framework -- in~\cite{pablo2025dd29}, HPS has been identified with a specific multigrid V-cycle.

Our primary goal below is to separate the \emph{discretization method} and the way in which the method treats the \emph{corner points}, thereby helping to build the modular view of HPS. For that reason we choose the standard Poisson problem on a rectangle and use the Q1-FEM discretization on a tensor product grid, where the basis functions firmly couple the corner point DoFs with others. In the HPS community, this would be considered a major issue as it prevents a straightforward definition of the local DtN. However, having discretized we show this can be resolved with little extra effort. We are not aware of the HPS and Q1-FEM coupling (or other simple low-order FEM) anywhere in the literature; this set-up should also provide a simple entry point into HPS methods for broader audience and FEM enthusiasts will notice that we do not rely on the Q1 elements in any way. We again highlight that in \emph{different communities} and in different context similar ideas already exists, see, e.g.~\cite{lu2024compression}, where the authors consider mixed-order curl-conforming FEM discretization for time-harmonic Maxwell equations in $\mathbb{R}^3$ -- an involved setting in which HPS is not mentioned but corner points and edge points are considered.

\section{The HPS method with corners}\label{sec_Q1HPS}

As noted above, we consider the simplest model problem
\begin{equation}\label{eqn_secHPS_PoissPrblm}
  \Delta u =f \quad \mathrm{in} \; \Omega:=(\alpha,\beta)\times(\gamma,\delta) \quad \mathrm{and} \quad
  u=g \quad \mathrm{on} \; \partial \Omega ,
\end{equation}
\noindent and start by outlining the structure of a general HPS method:

\begin{enumerate}
\item \emph{Partition} -- partition the domain $\Omega$ into subdomains.
\item \emph{Discretization \& Assembly} -- formulate, discretize and assemble the subdomain solution boundary operators for the subdomains.
\item \emph{Merge} --  merge the neighboring solution boundary operators and store the result.
\item \emph{Recursion} --  recurse and continue merging until we reach the entire domain $\Omega$.
\item \emph{Application} -- given data, apply the global solution boundary operator and calculate the solution on the boundaries of the subdomains.
\item \emph{Reconstruction} -- reconstruct the solution in the subdomains from the boundaries.
\end{enumerate}

As per the \emph{partition} stage, we the standard, grid-like set-up
\begin{equation*}
  \Omega_{e} = [a^{(x_1)}_{e}, b^{(x_1)}_{e}] \times [a^{(x_2)}_{e}, b^{(x_2)}_{e}] \subset \Omega,
\end{equation*}
\noindent see Figure~\ref{fig_secAssmbly_OmgPartition}-right, forming a non-overlapping decomposition of $\Omega$ with corner points; other decompositions can be treated identically~\cite[Figure 2]{lu2024compression}.

\subsection{The \emph{discretization \& assembly} stage}

\textbf{The analytical background.} We are interested in constructing the \emph{subdomain solution boundary operators} -- dealing with the Poisson problem, those are the subdomain DtNs. Let $u_{e}$ denote the solution on the subdomain $\Omega_{e}$, i.e.,
\begin{equation}\label{eqn_secAssmbly_LocalPoissonPrblm}
-\Delta u_{e} = f  \quad \mathrm{in}\; \Omega_e, \quad \mathrm{and} \quad
u_{e} = g_e  \quad \mathrm{on}\; \partial\Omega_e.
\end{equation}
\noindent We can split $u_{e}$ into the sum of the harmonic lift of the boundary data $g_e$, denoted by $u^{(g)}_e$ and the particular solution of the interior load $f$, denoted by $u^{(f)}_e$, obtaining
\begin{equation*}
\resizebox{.97\textwidth}{!}{$
-\Delta u^{(g)}_e = 0 \; \mathrm{in} \; \Omega_{e} \; \& \; u^{(g)}_e = g_e \; \mathrm{on} \; \partial\Omega_{e}
\quad \mathrm{and} \quad
-\Delta u^{(f)}_e = f \; \mathrm{in} \; \Omega_{e} \; \& \; u^{(f)}_e = 0 \; \mathrm{on}\; \partial\Omega_{e}.
$}
\end{equation*}
\noindent Analogously, we also split the Neumann trace of the solution, denoted by $\partial_n u_{e}$,
\begin{equation}\label{eqn_secAssmbly_LocalFlxsAsDtN_ContFromultn}
\partial_n u_{e} = \partial_n u^{(g)}_e + \partial_n u^{(f)}_e =:
\Lambda_{e} g_e + q_{e}
\end{equation}
\noindent featuring homogeneous DtN $\Lambda_{e}$ and the particular Neumann trace $q_{e}$ on $\Omega_{e}$.

\begin{figure}[t]
\centering
\resizebox{.99\textwidth}{!}{
\begin{tikzpicture}[thick] 
\draw (-3,0) -- (3,0) -- (3,2) -- (-3,2) -- (-3,0);
\draw (0,0) -- (0,2);
\draw (-1.5,2.45) node(one) {$\Omega_1$};
\draw (1.5,2.45) node(one) {$\Omega_2$};
\draw[dashed, opacity=0.55] (-3,.5) -- (3,.5);
\draw[dashed, opacity=0.55] (-3,1) -- (3,1);
\draw[dashed, opacity=0.55] (-3,1.5) -- (3,1.5);
\draw[dashed, opacity=0.55] (-2.5,0) -- (-2.5,2);
\draw[dashed, opacity=0.55] (-2,0) -- (-2,2);
\draw[dashed, opacity=0.55] (-1.5,0) -- (-1.5,2);
\draw[dashed, opacity=0.55] (-1,0) -- (-1,2);
\draw[dashed, opacity=0.55] (-.5,0) -- (-.5,2);
\draw[dashed, opacity=0.55] (0,0) -- (0,2);
\draw[dashed, opacity=0.55] (.5,0) -- (.5,2);
\draw[dashed, opacity=0.55] (1,0) -- (1,2);
\draw[dashed, opacity=0.55] (1.5,0) -- (1.5,2);
\draw[dashed, opacity=0.55] (2,0) -- (2,2);
\draw[dashed, opacity=0.55] (2.5,0) -- (2.5,2);
\draw[rounded corners, pattern={Lines[angle=45,distance={3pt/0.6}]},pattern color=blue, opacity=0.65] (-3.15, 0.35) rectangle (-2.85, 1.65);
\draw[rounded corners] (-3.15, 0.35) rectangle (-2.85, 1.65);
\draw[rounded corners, pattern={Lines[angle=45,distance={3pt/1.6}]},pattern color=blue, opacity=0.65] (-0.15, 0.35) rectangle (0.15, 1.65);
\draw[rounded corners] (-0.15, 0.35) rectangle (0.15, 1.65);
\draw[rounded corners, pattern={Lines[angle=135,distance={3pt/1.3}]},pattern color=green, opacity=0.65] (-0.15, 0.35) rectangle (0.15, 1.65);
\draw[rounded corners] (-0.15, 0.35) rectangle (0.15, 1.65);
\draw[rounded corners, pattern={Lines[angle=135,distance={3pt/0.6}]},pattern color=green, opacity=0.65] (2.85, 0.35) rectangle (3.15, 1.65);
\draw[rounded corners] (2.85, 0.35) rectangle (3.15, 1.65);
\draw[rounded corners, pattern={Lines[angle=45,distance={3pt/0.6}]},pattern color=blue, opacity=0.65] (-2.65, -0.15) rectangle (-0.35, 0.15);
\draw[rounded corners] (-2.65, -0.15) rectangle (-0.35, 0.15);
\draw[rounded corners, pattern={Lines[angle=45,distance={3pt/0.6}]},pattern color=blue, opacity=0.65] (-2.65, 1.85) rectangle (-0.35, 2.15);
\draw[rounded corners] (-2.65, 1.85) rectangle (-0.35, 2.15);
\draw[rounded corners, pattern={Lines[angle=135,distance={3pt/0.6}]},pattern color=green, opacity=0.65] (0.35, -0.15) rectangle (2.65, 0.15);
\draw[rounded corners] (0.35, -0.15) rectangle (2.65, 0.15);
\draw[rounded corners, pattern={Lines[angle=135,distance={3pt/0.6}]},pattern color=green, opacity=0.65] (0.35, 1.85) rectangle (2.65, 2.15);
\draw[rounded corners] (0.35, 1.85) rectangle (2.65, 2.15);
\draw[pattern={Lines[angle=45,distance={3pt/0.6}]},pattern color=blue, opacity=0.65] (-2.6, 0.4) rectangle (-0.4, 1.6);
\draw[] (-2.6, 0.4) rectangle (-0.4, 1.6);
\draw[pattern={Lines[angle=135,distance={3pt/0.6}]},pattern color=green, opacity=0.65] (0.4, 0.4) rectangle (2.6, 1.6);
\draw[] (0.4, 0.4) rectangle (2.6, 1.6);
\draw (-3.35,0.75) node(one) {{\scriptsize $\iota_1^{L}$}};
\draw (-2.45,2.4) node(one) {{\scriptsize $\iota_1^{T}$}};
\draw (-2.45,-0.4) node(one) {{\scriptsize $\iota_1^{B}$}};
\draw[rounded corners,draw=none,fill=white,fill opacity=0.99] (-1.95,1.0) rectangle (-1.55,1.5);
\draw (-1.75,1.25) node(one) {{\footnotesize $\iota_1^{\mathrm{int}}$}};
\draw (3.35,0.75) node(one) {{\scriptsize $\iota_2^{R}$}};
\draw (2.45,2.4) node(one) {{\scriptsize $\iota_2^{T}$}};
\draw (2.45,-0.4) node(one) {{\scriptsize $\iota_2^{B}$}};
\draw[rounded corners,draw=none,fill=white,fill opacity=0.99] (1.55,1.0) rectangle (1.95,1.5);
\draw (1.75,1.25) node(one) {{\footnotesize $\iota_2^{\mathrm{int}}$}};
\draw[rounded corners,draw=none,fill=white,fill opacity=0.99] (-0.135,0.45) rectangle (0.135,0.95);
\draw (0,0.75) node(one) {{\footnotesize $\iota_{1}^{I}$}};
\draw[rounded corners,draw=none,fill=white,fill opacity=0.99] (-0.135,1.05) rectangle (0.135,1.45);
\draw (0,1.25) node(one) {{\footnotesize $\iota_{2}^{I}$}};
\draw (-3,0) node(x) {{\tiny $\times$}};
\draw (-3,.5) node(x) {{\tiny $\times$}};
\draw (-3,1) node(x) {{\tiny $\times$}};
\draw (-3,1.5) node(x) {{\tiny $\times$}};
\draw (-3,2) node(x) {{\tiny $\times$}};
\draw (-2.5,0) node(x) {{\tiny $\times$}};
\draw (-2.5,.5) node(x) {{\tiny $\times$}};
\draw (-2.5,1) node(x) {{\tiny $\times$}};
\draw (-2.5,1.5) node(x) {{\tiny $\times$}};
\draw (-2.5,2) node(x) {{\tiny $\times$}};
\draw (-2,0) node(x) {{\tiny $\times$}};
\draw (-2,.5) node(x) {{\tiny $\times$}};
\draw (-2,1) node(x) {{\tiny $\times$}};
\draw (-2,1.5) node(x) {{\tiny $\times$}};
\draw (-2,2) node(x) {{\tiny $\times$}};
\draw (-1.5,0) node(x) {{\tiny $\times$}};
\draw (-1.5,.5) node(x) {{\tiny $\times$}};
\draw (-1.5,1) node(x) {{\tiny $\times$}};
\draw (-1.5,1.5) node(x) {{\tiny $\times$}};
\draw (-1.5,2) node(x) {{\tiny $\times$}};
\draw (-1,0) node(x) {{\tiny $\times$}};
\draw (-1,.5) node(x) {{\tiny $\times$}};
\draw (-1,1) node(x) {{\tiny $\times$}};
\draw (-1,1.5) node(x) {{\tiny $\times$}};
\draw (-1,2) node(x) {{\tiny $\times$}};
\draw (-.5,0) node(x) {{\tiny $\times$}};
\draw (-.5,.5) node(x) {{\tiny $\times$}};
\draw (-.5,1) node(x) {{\tiny $\times$}};
\draw (-.5,1.5) node(x) {{\tiny $\times$}};
\draw (-.5,2) node(x) {{\tiny $\times$}};
\draw (0,0) node(x) {{\tiny $\times$}};
\draw (0,.5) node(x) {{\tiny $\times$}};
\draw (0,1) node(x) {{\tiny $\times$}};
\draw (0,1.5) node(x) {{\tiny $\times$}};
\draw (0,2) node(x) {{\tiny $\times$}};
\draw (.5,0) node(x) {{\tiny $\times$}};
\draw (.5,.5) node(x) {{\tiny $\times$}};
\draw (.5,1) node(x) {{\tiny $\times$}};
\draw (.5,1.5) node(x) {{\tiny $\times$}};
\draw (.5,2) node(x) {{\tiny $\times$}};
\draw (1,0) node(x) {{\tiny $\times$}};
\draw (1,.5) node(x) {{\tiny $\times$}};
\draw (1,1) node(x) {{\tiny $\times$}};
\draw (1,1.5) node(x) {{\tiny $\times$}};
\draw (1,2) node(x) {{\tiny $\times$}};
\draw (1.5,0) node(x) {{\tiny $\times$}};
\draw (1.5,.5) node(x) {{\tiny $\times$}};
\draw (1.5,1) node(x) {{\tiny $\times$}};
\draw (1.5,1.5) node(x) {{\tiny $\times$}};
\draw (1.5,2) node(x) {{\tiny $\times$}};
\draw (2,0) node(x) {{\tiny $\times$}};
\draw (2,.5) node(x) {{\tiny $\times$}};
\draw (2,1) node(x) {{\tiny $\times$}};
\draw (2,1.5) node(x) {{\tiny $\times$}};
\draw (2,2) node(x) {{\tiny $\times$}};
\draw (2.5,0) node(x) {{\tiny $\times$}};
\draw (2.5,.5) node(x) {{\tiny $\times$}};
\draw (2.5,1) node(x) {{\tiny $\times$}};
\draw (2.5,1.5) node(x) {{\tiny $\times$}};
\draw (2.5,2) node(x) {{\tiny $\times$}};
\draw (3,0) node(x) {{\tiny $\times$}};
\draw (3,.5) node(x) {{\tiny $\times$}};
\draw (3,1) node(x) {{\tiny $\times$}};
\draw (3,1.5) node(x) {{\tiny $\times$}};
\draw (3,2) node(x) {{\tiny $\times$}};
\draw (4,0) -- (8,0) -- (8,2) -- (4,2) -- (4,0);
\draw (5,0) -- (5,2);
\draw (6,0) -- (6,2);
\draw (7,0) -- (7,2);
\draw (4,0.5) -- (8,0.5);
\draw (4,1) -- (8,1);
\draw (4,1.5) -- (8,1.5);
\draw (6,2.25) node(one) {$\Omega$};
\draw (4.55,0.21) node(one) {{\footnotesize $\Omega_{1,1}$}};
\draw (5.55,0.21) node(one) {{\footnotesize $\Omega_{1,2}$}};
\draw (6.55,0.21) node(one) {{\footnotesize $\Omega_{1,3}$}};
\draw (7.55,0.21) node(one) {{\footnotesize $\Omega_{1,4}$}};
\draw (4.55,0.77) node(one) {{\footnotesize $\Omega_{2,1}$}};
\draw (5.55,0.76) node(one) {{\footnotesize $\Omega_{2,2}$}};
\draw (6.55,0.72) node(one) {{\footnotesize $\cdots$}};
\draw[fill=gray!10] (5,0.25) circle (0.15cm);
\node at (5,0.25) {{\tiny $1$}};
\draw[fill=gray!10] (7,0.25) circle (0.15cm);
\node at (7,0.25) {{\tiny $2$}};
\draw[fill=gray!10] (5,0.75) circle (0.15cm);
\node at (5,0.75) {{\tiny $3$}};
\draw[fill=gray!10] (7,0.75) circle (0.15cm);
\node at (7,0.75) {{\tiny $4$}};
\draw[fill=gray!10] (5,1.25) circle (0.15cm);
\node at (5,1.25) {{\tiny $5$}};
\draw[fill=gray!10] (7,1.25) circle (0.15cm);
\node at (7,1.25) {{\tiny $6$}};
\draw[fill=gray!10] (5,1.75) circle (0.15cm);
\node at (5,1.75) {{\tiny $7$}};
\draw[fill=gray!10] (7,1.75) circle (0.15cm);
\node at (7,1.75) {{\tiny $8$}};
\draw[fill=gray!10] (4.25,0.4) rectangle (4.75,0.6);
\node at (4.5,0.5) {{\tiny $9$}};
\draw[fill=gray!10] (5.25,0.4) rectangle (5.75,0.6);
\node at (5.5,0.5) {{\tiny $10$}};
\draw[fill=gray!10] (4.25,1.4) rectangle (4.75,1.6);
\node at (4.5,1.5) {{\tiny $11$}};
\draw[fill=gray!10] (5.25,1.4) rectangle (5.75,1.6);
\node at (5.5,1.5) {{\tiny $12$}};
\draw[fill=gray!10] (6.25,0.4) rectangle (6.75,0.6);
\node at (6.5,0.5) {{\tiny $13$}};
\draw[fill=gray!10] (7.25,0.4) rectangle (7.75,0.6);
\node at (7.5,0.5) {{\tiny $14$}};
\draw[fill=gray!10] (6.25,1.4) rectangle (6.75,1.6);
\node at (6.5,1.5) {{\tiny $15$}};
\draw[fill=gray!10] (7.25,1.4) rectangle (7.75,1.6);
\node at (7.5,1.5) {{\tiny $16$}};
\end{tikzpicture}
}
\caption{Left: two neighbouring subdomains $\Omega_1,\Omega_2$ with the index sets of the grids. We have $\iota_1^I\equiv \iota_1^R$ and $\iota_2^I\equiv \iota_2^L$ and also see the corner index sets $\iota_1^C,\iota_2^C$, although not separately highlighted. Finally, we set $\iota^{\partial}_{e} := 
\iota^{L}_{e} \cup \iota^{R}_{e} \cup \iota^{T}_{e} \cup \iota^{B}_{e} \cup \iota^{C}_{e}$ so that $\iota_{e} = \iota^{\mathrm{int}}_{e} \cup \iota^{\partial}_{e}$ for any $e$. Right: (incomplete) illustration of the nested dissection merge hierarchy ordering, see~\cite[Appendix A]{Martinsson2013}.}\label{fig_secAssmbly_OmgPartition}
\end{figure}
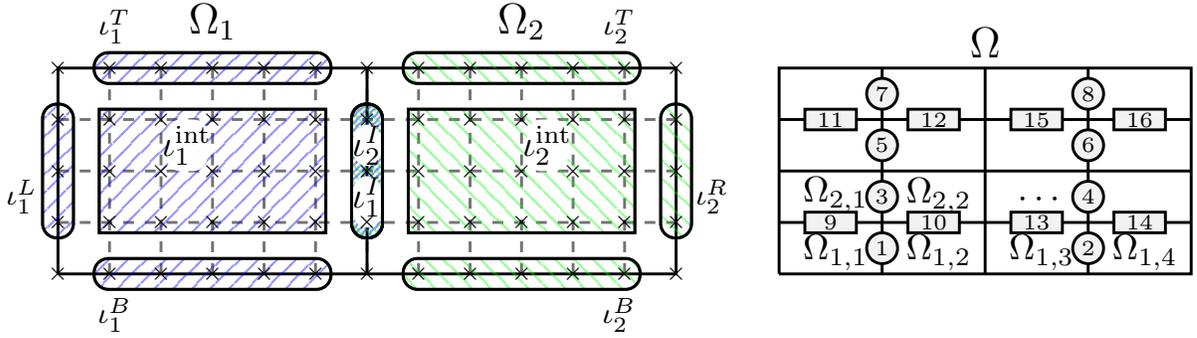

\textbf{Discretization.} We first introduce grid nodes in $\Omega^{e}$ in a tensor-product manner along the $x_1$ and $x_2$ axis. On this grid we consider the $Q_1$ finite element discretization of~\eqref{eqn_secAssmbly_LocalPoissonPrblm} and index the local DoFs by integer pairs $\iota_{e}=(i,j)$; see Figure~\ref{fig_secAssmbly_OmgPartition} for the details. Applying integration by parts to the continuous weak form of~\eqref{eqn_secAssmbly_LocalPoissonPrblm} gives
\begin{equation*}
  \int_{\Omega_e} \nabla u_{e} \cdot \nabla \phi_m\,\mathrm{d}\mathbf{x}
  = \int_{\Omega_e} f\,\phi_m\,\mathrm{d}\mathbf{x}
  + \int_{\partial\Omega_e} (\partial_n u_{e})\,\phi_m\,\mathrm{d}s, \quad m\in \iota_{e}
\end{equation*}
\noindent and then, after approximating $u_{e}$ in the Q1-FEM basis and reordering the DoFs, we get the discretized system for the unknown coefficients $\mathbf{u}_e$
\begin{equation}\label{eqn_secAssmbly_DiscrtzdIntgrtnByParts}
\begin{bmatrix}
    A_{e} & B_{e}\\[3pt]
    C_{e} & D_{e}
\end{bmatrix}
\begin{bmatrix}
    \mathbf{u}_{e}^{\mathrm{int}}\\[3pt]
    \mathbf{u}_{e}^{\partial}
\end{bmatrix}
  =
\begin{bmatrix}
    \mathbf{f}_{e}^{\mathrm{int}}\\[3pt]
    \mathbf{f}_{e}^{\partial}
\end{bmatrix} +
\begin{bmatrix}
	\mathbf{0} \\[3pt]
	\partial_n \mathbf{u}_{e}
\end{bmatrix},
\end{equation}
\noindent where $\partial_n \mathbf{u}_{e}$ is the Q1-FEM discretization of the Neumann trace of $\partial_n u_{e}$ at $\iota_{e}^{\mathrm{int}}$. 

\textbf{Assembly.} Equation~\eqref{eqn_secAssmbly_DiscrtzdIntgrtnByParts} shows that the \emph{negative residual along the interface} is the finite element representation of the approximate normal fluxes along the boundary. As these fluxes are unknown, the equations in the second block-row give the formula,
\begin{equation*}
\partial_n \mathbf{u}_{e} = C_e \mathbf{u}_{e}^{\mathrm{int}} + D_e \mathbf{u}_{e}^{\partial} - \mathbf{f}_{e}^{\partial},
\end{equation*}
\noindent which, after elimination of the interior DoFs, becomes
\begin{equation}\label{eqn_secAssmbly_LocalFlxsAsDtN_DiscrtMirroring}
\partial_n \mathbf{u}_{e} \equiv \mathbf{r}_{e}^{\partial} = \left( D_{e} - C_{e} A_{e}^{-1} B_{e} \right) \mathbf{u}_{e}^{\partial} + C_{e}A_{e}^{-1} \mathbf{f}_{e}^{\mathrm{int}} - \mathbf{f}_{e}^{\partial} =: S_{e} \mathbf{u}_{e}^{\partial} + \mathbf{h}_{e}.
\end{equation}
\noindent This relation mirrors the continuous decomposition~\eqref{eqn_secAssmbly_LocalFlxsAsDtN_ContFromultn}:
the Schur complement $S_{e}$ acts as the discrete homogeneous DtN operator mapping boundary values $\mathbf{u}_{e}^{\partial}$ to their induced boundary fluxes, while $\mathbf{h}_{e}$ represents the discrete flux produced by the interior load under homogeneous Dirichlet conditions. The \emph{assembly} stage of HPS methods consists of computing (or approximating) the matrices $S_{e}$ so that the right-hand side of~\eqref{eqn_secAssmbly_LocalFlxsAsDtN_DiscrtMirroring} can be evaluated \emph{rapidly} in the \emph{application} stage.

\subsection{The \emph{merge} stage}\label{sec_MergeStep}
Having two subdomains, say $\Omega_1,\Omega_2$, with finished \emph{assembly} stage that share an interface, we want to assemble the solution boundary operator for $\Omega_1 \cup\Omega_2$. The ordering in which we will pick the subdomain pairs matters as it highly influences the parallelizability of the resulting solver; we follow the nested dissection ordering as illustrated in Figure~\ref{fig_secAssmbly_OmgPartition}-right; first we merge \emph{horizontally} and then \emph{vertically}.

\textbf{Horizontal merge (left $\Omega_1$ and right $\Omega_2$).} The true solution of~\eqref{eqn_secHPS_PoissPrblm} is continuous and has balanced fluxes (i.e., residuals) across the interface, i.e.,
\begin{equation}\label{eqn_secMerge_Hor_ContOfSolBalancOfFlxs}
\mathbf{u}_1(\iota_1^{I}) = \mathbf{u}_2(\iota_2^{I}) \quad \mathrm{and} \quad \mathbf{r}_1(\iota_1^{I}) + \mathbf{r}_2(\iota_2^{I}) = 0.
\end{equation}
\noindent Recalling~\eqref{eqn_secAssmbly_LocalFlxsAsDtN_DiscrtMirroring} and blocking it according to $\iota_e^{\partial} = \iota_e^{I} \cup \iota_e^{\partial\backslash I}$ for $e=1,2$ gives
\begin{equation}\label{eqn_secMerge_Hor_GnrlResdlEqnAlongIntrfc}
\begin{bmatrix}
\mathbf{r}_e(\iota_e^{\partial\setminus I})\\
\mathbf{r}_e(\iota_e^{I})
\end{bmatrix}
=
\begin{bmatrix}
S_e(\iota_e^{\partial\setminus I},\iota_e^{\partial\setminus I}) & S_e(\iota_e^{\partial\setminus I},\iota_e^{I})\\
S_e(\iota_e^{I},\iota_e^{\partial\setminus I}) & S_e(\iota_e^{I},\iota_e^{I})
\end{bmatrix}
\begin{bmatrix}
\mathbf{u}_e(\iota_e^{\partial\setminus I})\\
\mathbf{u}_e(\iota_e^{I})
\end{bmatrix}
+
\begin{bmatrix}
\mathbf{h}_e(\iota_e^{\partial\setminus I})\\
\mathbf{h}_e(\iota_e^{I})
\end{bmatrix}.
\end{equation}
\noindent Summing the second block-rows for $e=1,2$, using~\eqref{eqn_secMerge_Hor_ContOfSolBalancOfFlxs} and reordering gives
\begin{equation}\label{eqn_secMerge_Hor_SolOnIntrfcIndsFrml}
\left( S_1(\iota_1^{I},\iota_1^{I}) + S_2(\iota_2^{I},\iota_2^{I}) \right) \mathbf{u}_1(\iota_1^{I}) =
- \sum\limits_{e=1,2} \mathbf{h}_e(\iota_e^{I}) + S_e(\iota_e^{I},\iota_e^{\partial\setminus I})\mathbf{u}_e(\iota_e^{\partial\setminus I}).
\end{equation}
\noindent Returning to~\eqref{eqn_secMerge_Hor_GnrlResdlEqnAlongIntrfc}, we concatenate the equations for the residuals on the ``merged boundary'' $r_e(\iota_e^{\partial\setminus I}), e=1,2$, use~\eqref{eqn_secMerge_Hor_SolOnIntrfcIndsFrml} in both and reorder so as to obtain
\begin{equation}\label{eqn_secMerge_Hor_ResOnExteriorIndsAsDtN}
\begin{bmatrix}
\mathbf{r}_1(\iota_1^{\partial\setminus I})\\[2pt]
\mathbf{r}_2(\iota_2^{\partial\setminus I})
\end{bmatrix}
=
S^{\mathrm{H}}
\begin{bmatrix}
\mathbf{u}_1(\iota_1^{\partial\setminus I})\\[2pt]
\mathbf{u}_2(\iota_2^{\partial\setminus I})
\end{bmatrix}
+
\mathbf{h}^{\mathrm{H}},
\end{equation}
\noindent i.e., the \emph{horizontally} merged boundary solution operators as in~\eqref{eqn_secAssmbly_LocalFlxsAsDtN_DiscrtMirroring} with
\begin{equation*}
\begin{aligned}
S^{\mathrm{H}}
=&
\begin{bmatrix}
S_1(\iota_1^{\partial\setminus I},\iota_1^{\partial\setminus I}) & 0\\
0 & S_2(\iota_2^{\partial\setminus I},\iota_2^{\partial\setminus I})
\end{bmatrix}
- \\
&\begin{bmatrix}
S_1(\iota_1^{\partial\setminus I},\iota_1^{I})\\
S_2(\iota_2^{\partial\setminus I},\iota_2^{I})
\end{bmatrix}
\big(S_1(\iota_1^{I},\iota_1^{I}) + S_2(\iota_2^{I},\iota_2^{I})\big)^{-1}
\begin{bmatrix}
S_1(\iota_1^{I},\iota_1^{\partial\setminus I}) & S_2(\iota_2^{I},\iota_2^{\partial\setminus I})
\end{bmatrix},\\[1ex]
\mathbf{h}^{\mathrm{H}}
=&
\begin{bmatrix}
\mathbf{h}_1(\iota_1^{\partial\setminus I})\\[2pt]
\mathbf{h}_2(\iota_2^{\partial\setminus I})
\end{bmatrix} -
\begin{bmatrix}
S_1(\iota_1^{\partial\setminus I},\iota_1^{I})\\[2pt]
S_2(\iota_2^{\partial\setminus I},\iota_2^{I})
\end{bmatrix}
\big(S_1(\iota_1^{I},\iota_1^{I}) + S_2(\iota_2^{I},\iota_2^{I})\big)^{-1}
\big(\mathbf{h}_1(\iota_1^{I}) + \mathbf{h}_2(\iota_2^{I})\big).
\end{aligned}
\end{equation*}

\textbf{Vertical merge (bottom $\Omega_1$ and top $\Omega_2$).} Say we have ``horizontally'' merged the boundary solution operators for two couples of subdomains $\Omega_{1L},\Omega_{1R}$ and $\Omega_{2L}, \Omega_{2R}$, e.g. the merges \textcircled{1} and \textcircled{3} in Figure~\ref{fig_secAssmbly_OmgPartition}, and we are ready to merge along the vertical interface -- labeled \fbox{{\scriptsize 9}} and \fbox{{\scriptsize 10}} -- \emph{and then also at the corner point enclosed between the already merged interfaces}. First, keeping the enclosed corner DoF, indexed\footnote{We use $\iota^{C}_{\Gamma}\equiv c$ as an \emph{absolute} index across the indexing in the four subdomains $\Omega_{1L}, \Omega_{1R}, \Omega_{2L}, \Omega_{2R}$, even though the point has likely a different index in each of them.} by $\iota^{C}_{\Gamma}\equiv c$, uneliminated, the steps in merges \fbox{{\scriptsize 9}} and \fbox{{\scriptsize 10}} carry through identically to the horizontal merges \textcircled{1} or \textcircled{3}, only now the index set $\iota_e^{I}$ is \emph{disjoint}, e.g., $\iota_1^{I} = \iota_{1L}^{T} \cup \iota_{1R}^{T}$, and the index sets $\iota_1^{\partial\backslash I},\iota_2^{\partial\backslash I}$ contain the enclosed corner index $\iota^{C}_{\Gamma}$. That is, we have
\begin{equation}\label{eqn_secMerge_Ver_ResOnExteriorIndsAsDtN}
\begin{bmatrix}
\mathbf{r}_1(\iota_1^{\partial\setminus I})\\[2pt]
\mathbf{r}_2(\iota_2^{\partial\setminus I})
\end{bmatrix}
=
S^{\mathrm{V}}
\begin{bmatrix}
\mathbf{u}_1(\iota_1^{\partial\setminus I})\\[2pt]
\mathbf{u}_2(\iota_2^{\partial\setminus I})
\end{bmatrix}
+
\mathbf{h}^{\mathrm{V}},
\end{equation}
\noindent with the \emph{vertically} merged boundary solution operators as in~\eqref{eqn_secAssmbly_LocalFlxsAsDtN_DiscrtMirroring}, i.e.,
\begin{equation*}
\begin{aligned}
S^{\mathrm{V}}
=&
\begin{bmatrix}
S_1(\iota_1^{\partial\setminus I},\iota_1^{\partial\setminus I}) & 0\\
0 & S_2(\iota_2^{\partial\setminus I},\iota_2^{\partial\setminus I})
\end{bmatrix}
- \\
&\begin{bmatrix}
S_1(\iota_1^{\partial\setminus I},\iota_1^{I})\\
S_2(\iota_2^{\partial\setminus I},\iota_2^{I})
\end{bmatrix}
\big(S_1(\iota_1^{I},\iota_1^{I}) + S_2(\iota_2^{I},\iota_2^{I})\big)^{-1}
\begin{bmatrix}
S_1(\iota_1^{I},\iota_1^{\partial\setminus I}) & S_2(\iota_2^{I},\iota_2^{\partial\setminus I})
\end{bmatrix},\\[1ex]
\mathbf{h}^{\mathrm{V}}
=&
\begin{bmatrix}
\mathbf{h}_1(\iota_1^{\partial\setminus I})\\[2pt]
\mathbf{h}_2(\iota_2^{\partial\setminus I})
\end{bmatrix} -
\begin{bmatrix}
S_1(\iota_1^{\partial\setminus I},\iota_1^{I})\\[2pt]
S_2(\iota_2^{\partial\setminus I},\iota_2^{I})
\end{bmatrix}
\big(S_1(\iota_1^{I},\iota_1^{I}) + S_2(\iota_2^{I},\iota_2^{I})\big)^{-1}
\big(\mathbf{h}_1(\iota_1^{I}) + \mathbf{h}_2(\iota_2^{I})\big).
\end{aligned}
\end{equation*}

\textbf{Corner merge (corner of $\Omega_{1L}, \Omega_{1R}, \Omega_{2L}, \Omega_{2R}$).} Analogously to~\eqref{eqn_secMerge_Hor_ContOfSolBalancOfFlxs}, we have
\begin{equation}\label{eqn_secMerge_Cor_ContOfSolAndBalancOfFlxs}
\mathbf{u}_{1}(c)=\mathbf{u}_{2}(c)=\mathbf{u}(c)
\quad \mathrm{and} \quad
\sum\limits_{e=1,2} \mathbf{r}_e(c) = 0.
\end{equation}
\noindent Collecting the $\iota^{C}_{\Gamma}\equiv c$ equations from~\eqref{eqn_secMerge_Ver_ResOnExteriorIndsAsDtN} and inserting them into~\eqref{eqn_secMerge_Cor_ContOfSolAndBalancOfFlxs} gives
\begin{equation*}
\sum_{e=1,2} S_e(c,c)\,\mathbf{u}(c) =
- \sum_{e=1,2} S_e(c,E_e)\mathbf{u}_e(E_e) + \mathbf{h}_e(c),
\end{equation*}
\noindent where $E_e$ are the indices of points on $\partial\Omega_e \cap \partial(\Omega_1\cup \Omega_2)$. Solving for $\mathbf{u}(c)$ and inserting back in~\eqref{eqn_secMerge_Ver_ResOnExteriorIndsAsDtN} gives a system on the exterior index set $E:= E_1\cup E_2$ for fluxes
\begin{equation*}
\begin{gathered}
\mathbf{r}(E) =
S^{\mathrm{corner}}\,\mathbf{u}(E) + \mathbf{h}^{\mathrm{corner}} \; \mathrm{with}\\
\begin{aligned}
S^{\mathrm{corner}} &=
S^{\mathrm{V}}(E,E) -
S^{\mathrm{V}}(E,c) \Big(S^{\mathrm{V}}(c,c)\Big)^{-1} S^{\mathrm{V}}(c,E), \\
\mathbf{h}^{\mathrm{corner}} &=
\mathbf{h}^{\mathrm{V}}(E) - S^{\mathrm{V}}(E,c) \Big(S^{\mathrm{V}}(c,c)\Big)^{-1} \mathbf{h}^{\mathrm{V}}(c),
\end{aligned}
\end{gathered}
\end{equation*}
\noindent with identical structure of the resulting boundary solution operator as in~\eqref{eqn_secMerge_Ver_ResOnExteriorIndsAsDtN} or~\eqref{eqn_secMerge_Hor_ResOnExteriorIndsAsDtN}.

\textbf{What is the point?} 
First, this treatment of the corner points is fundamentally different to the ``change of basis'' approach used, e.g., in~\cite{Martinsson2015} -- no retabulation, rather following the same ground ideas behind HPS. Second, it is also fundamentally different from the ``ignore'' approach used, e.g., in~\cite{GeldermansGillman2019}, as that is simply not an option due to the corner point DoFs coupling. Third, this is in fact very similar to~\cite{Martinsson2013}, but outlined only within the SEM context with the aforementioned benefits. Our point is that the corner points should be merged once the surrounding interfaces have been merged to maximize the efficiency and doing that follows analogous steps used before, even when using fully coupled corner DoFs of FEM.

\subsection{The \emph{recursion} stage}
Having successfully eliminated the interface and enclosed corner DoFs, we recurse and continue until reaching a problem on $\partial \Omega$, where the Dirichlet trace is known. Let $\Omega = \bigcup_e \Omega_e$ be the decomposition into subdomains, each with its local DtN $(S_e,\mathbf{h}_e)$. The global system on the \emph{skeleton} -- i.e., on $\bigcup_e \partial\Omega_e$ -- reads
\begin{equation*}
  S\, \mathbf{u}^{(\partial)} = \mathbf{h}^{(\partial)},
\end{equation*}
\noindent where $S$ is built from $S_e$ based on the interface continuity and flux balances. We order the boundary indices $\iota^{(\partial)} = \bigcup_e \iota_e^{\partial}$ by the merging hierarchy: first domain boundaries, then merged interfaces and enclosed corners (following the nested dissection ordering). Hence the \emph{to-internalize} indices come after the \emph{active} exterior indices,
\begin{equation*}
  \iota^{(\partial)} =
  \big(
    \iota^{(1)}_{\mathrm{ext}},\
    \iota^{(1)}_{\mathrm{merge}},\
    \iota^{(2)}_{\mathrm{merge}},\
    \dots,\
    \iota^{(L)}_{\mathrm{merge}}
  \big),
\end{equation*}
\noindent where $\iota^{(\ell)}_{\mathrm{merge}}$ are the indices eliminated at recursion level~$\ell$. Then $S$ has the structure
\begin{equation*}
  S =
  \begin{bmatrix}
    S_{EE} & S_{EM}\\
    S_{ME} & S_{MM}
  \end{bmatrix},
\end{equation*}
\noindent where $E$ and $M$ representing the exterior and to-be-merged blocks. The \emph{merge} step corresponds precisely to eliminating the $S_{MM}$ block via its Schur complement:
\begin{equation*}
  \widehat{S}_{EE} =
  S_{EE} - S_{EM}\,S_{MM}^{-1}\,S_{ME}
  \quad \mathrm{and} \quad 
  \widehat{\mathbf{h}}_E =
  \mathbf{h}_E - S_{EM}\,S_{MM}^{-1}\,\mathbf{h}_M,
\end{equation*}
\noindent where $(\widehat{S}_{EE},\widehat{\mathbf{h}}_E)$ defines the reduced DtN operator and right-hand side of the updated skeleton after that merge. Proceeding recursively the calculation always follows the same two-domain pattern, possibly extended by enclosed-corner junctions. That is the HPS skeleton solver can be interpreted as a single recursive Schur complement elimination applied to the global skeleton matrix $S$. Each recursion step in HPS corresponds to eliminating the block $(\iota^{(\ell)}_{\mathrm{merge}},\iota^{(\ell)}_{\mathrm{merge}})$ corresponding to indices merged at that level of the hierarchy. At the top of the recursion we get the final reduced operator $S^{(L)}$ on $\partial\Omega$, whose equilibrium equation represents the DtN map of $\Omega$.

\section{Numerical illustration}
\begin{table}[t]
\centering
\caption{Average speedup (left) and break-even solves (right) relative to MATLAB's backslash.}
\resizebox{1.\textwidth}{!}{%
\begin{tabular}{|c|cccc|cccc|}
  \multicolumn{1}{c}{}
  & \multicolumn{4}{c}{Speedup} & \multicolumn{4}{c}{Break-even solves} \\
  \hline
\diagbox[
  width=15em,
  height=3.5\baselineskip,
  font=\footnotesize
]{\# subdomains}{\makecell{\# elements p/subdomain}} &
$4\times4$ & $8\times8$ & $16\times16$ & $32\times32$ &
$4\times4$ & $8\times8$ & $16\times16$ & $32\times32$ \\
\hline
$4\times4$   & 1 & 3 & 5  & 12 & N/A & 4 & 2 & 2 \\
$8\times8$   & 2 & 6 & 9  & 17 & 5   & 4 & 3 & 2 \\
$16\times16$ & 3 & 7 & 21 & 28 & 5   & 5 & 4 & 3 \\
$32\times32$ & 4 & 11& 26 & 35 & 7   & 5 & 4 & 4 \\
\hline
\end{tabular}}
\end{table}
We conclude by showcasing the performance of HPS, implemented in MATLAB with very few optimizations. The build stage is computationally costly so the method is useful when we have several different right-hand sides. We run our tests on a laptop with 32GB RAM and a i7-12700H Intel microprocessor with six 4.7 GHz performance cores, eight 3.5GHz efficient cores and twenty total threads with performance-core hyperthreading. We take the MATLAB's backslash for the skeleton problem as the benchmark (ignoring the reconstruction). This comparison is stricter than a full solution comparison, where backslash would process a significantly larger operator.

\bibliographystyle{plainnat}
\bibliography{refs}

\begin{thebibliography}{26}
\providecommand{\natexlab}[1]{#1}
\providecommand{\url}[1]{\texttt{#1}}
\expandafter\ifx\csname urlstyle\endcsname\relax
  \providecommand{\doi}[1]{doi: #1}\else
  \providecommand{\doi}{doi: \begingroup \urlstyle{rm}\Url}\fi

\bibitem[Babb et~al.(2018)Babb, Gillman, Hao, and
  Martinsson]{BabbGillmanHaoMartinsson2018}
T.~Babb, A.~Gillman, S.~Hao, and P.-G. Martinsson.
\newblock An accelerated {P}oisson solver based on a multidomain spectral
  discretization.
\newblock \emph{{BIT} {N}um. {M}ath.}, 58\penalty0 (4):\penalty0 851--879,
  2018.

\bibitem[Beams et~al.(2020)Beams, Gillman, and Hewett]{BeamsGillmanHewett2020}
N.~N. Beams, A.~Gillman, and R.~J. Hewett.
\newblock A parallel shared-memory implementation of a high-order accurate
  solution technique for variable coefficient {H}elmholtz problems.
\newblock \emph{{C}omp. \& {M}ath. w. {A}ppl.}, 79\penalty0 (4):\penalty0
  996--1011, 2020.

\bibitem[B{\"{o}}rm et~al.(2006)B{\"{o}}rm, Grasedyck, and
  Hacksbusch]{borm2006hierarchicalmatrices}
S.~B{\"{o}}rm, L.~Grasedyck, and W.~Hacksbusch.
\newblock {H}ierarchical {M}atrices: Lecture notes no. 21, 2006.

\bibitem[Chen(2002)]{Chen2002}
Y.~Chen.
\newblock A fast, direct algorithm for the {L}ippmann–{S}chwinger integral
  equation in two dimensions.
\newblock \emph{{A}dv. in {C}omp. {M}ath.}, 16\penalty0 (2-3):\penalty0
  175--190, 2002.

\bibitem[Chipman et~al.(2024)Chipman, Calhoun, and
  Burstedde]{ChipmanCalhounBurstedde2024}
D.~Chipman, D.~Calhoun, and C.~Burstedde.
\newblock A fast direct solver for elliptic {PDE}s on a hierarchy of adaptively
  refined quadtrees.
\newblock \emph{arXiv:2402.14936}, 2024.

\bibitem[Claeys and Parolin(2022)]{ClaeysParolin2022}
X.~Claeys and E.~Parolin.
\newblock Robust treatment of cross-points in optimized schwartz methods.
\newblock \emph{{N}umer. {M}ath.}, 151:\penalty0 405--442, 2022.

\bibitem[Fortunato(2024)]{Fortunato2024}
D.~Fortunato.
\newblock A high-order fast direct solver for surface {PDE}s.
\newblock \emph{{SIAM} {J}. on {S}ci. {C}omp.}, 46\penalty0 (4):\penalty0
  A2582--A2606, 2024.

\bibitem[Gander and Santugini-Repiquet(2016)]{gander2014cross}
M.~J. Gander and K.~Santugini-Repiquet.
\newblock Cross-points in domain decomposition methods with a finite element
  discretization.
\newblock \emph{{ETNA}}, 45:\penalty0 219--240, 2016.

\bibitem[Geldermans and Gillman(2019)]{GeldermansGillman2019}
P.~Geldermans and A.~Gillman.
\newblock An adaptive high order direct solution technique for elliptic
  boundary value problems.
\newblock \emph{{SIAM} {J}. on {S}ci. {C}omp.}, 41\penalty0 (1):\penalty0
  A292--A315, 2019.

\bibitem[George(1973)]{george1973nested}
A.~George.
\newblock Nested dissection of a regular finite element mesh.
\newblock \emph{{SIAM} {J}. on {N}um. {A}nal.}, 10\penalty0 (2):\penalty0
  345--363, 1973.

\bibitem[Gillman(2011)]{Gillman2011}
A.~Gillman.
\newblock \emph{Fast direct solvers for elliptic partial differential
  equations}.
\newblock {Ph.D. Thesis}, University of {C}olorado {B}oulder, 2011.

\bibitem[Gillman and Martinsson(2014)]{GillmanMartinsson2014}
A.~Gillman and P.-G. Martinsson.
\newblock A direct solver with {$\mathcal{O}(N)$} complexity for variable
  coefficient elliptic {PDE}s discretized via a high-order composite spectral
  collocation method.
\newblock \emph{{SIAM} {J}. on {S}ci. {C}omp.}, 36\penalty0 (4):\penalty0
  A2023--A2046, 2014.

\bibitem[Gillman et~al.(2015)Gillman, Barnett, and
  Martinsson]{GillmanBarnettMartinsson2015}
A.~Gillman, A.~H. Barnett, and P.-G. Martinsson.
\newblock A spectrally accurate direct solution technique for frequency-domain
  scattering problems with variable media.
\newblock \emph{{BIT} {N}um. {M}ath.}, 55\penalty0 (1):\penalty0 141--170,
  2015.

\bibitem[Grasedyck et~al.(2009)Grasedyck, Kriemann, and
  Le~Borne]{GrasedyckKriemannLeBorne2009}
L.~Grasedyck, R.~Kriemann, and S.~Le~Borne.
\newblock Domain decomposition based {$\mathcal{H}$-LU} preconditioning.
\newblock \emph{{N}umerische {M}ath.}, 112\penalty0 (4):\penalty0 565--600,
  2009.

\bibitem[Greengard et~al.(2009)Greengard, Gueyffier, Martinsson, and
  Rokhlin]{GreengardGueyffierMartinssonRokhlin2009}
L.~Greengard, D.~Gueyffier, P.-G. Martinsson, and V.~Rokhlin.
\newblock Fast direct solvers for integral equations in complex
  three-dimensional domains.
\newblock \emph{{A}cta {N}um.}, 18:\penalty0 243–275, 2009.

\bibitem[Hackbusch(2015)]{hackbusch2015hierarchical}
W.~Hackbusch.
\newblock \emph{Hierarchical {M}atrices: {A}lgorithms and {A}nalysis}.
\newblock Springer {B}erlin, 2015.

\bibitem[Kump et~al.(2025)Kump, Yesypenko, and
  Martinsson]{KumpYesypenkoMartinsson2025}
Y.~Kump, A.~Yesypenko, and P.-G. Martinsson.
\newblock A two-level direct solver for the hierarchical
  {P}oincar{\'{e}}–{S}teklov method.
\newblock \emph{arXiv:2503.04033}, 2025.

\bibitem[Lu and Lee(2024)]{lu2024compression}
J.~Lu and J.-F. Lee.
\newblock A compression scheme for domain decomposition method in solving
  electromagnetic problems.
\newblock \emph{{J}ournal of {C}omputational {P}hysics}, 503:\penalty0 112824,
  2024.

\bibitem[Lucero~Lorca(2025)]{pablo2025dd29}
J.~P. Lucero~Lorca.
\newblock Towards a multigrid preconditioner interpretation of hierarchical
  {P}oincar{\'{e}}--steklov solvers.
\newblock \emph{Submitted in the proceedings of DD29}, 2025.

\bibitem[{Lucero Lorca} et~al.(2024){Lucero Lorca}, Beams, Beecroft, and
  Gillman]{LuceroLorcaBeamsBeecroftGillman2024}
J.~P. {Lucero Lorca}, N.~Beams, D.~Beecroft, and A.~Gillman.
\newblock An iterative solver for the {HPS} discretization applied to three
  dimensional {H}elmholtz problems.
\newblock \emph{{SIAM} {J}. {S}ci. {C}omp.}, 46\penalty0 (1):\penalty0
  A80--A104, 2024.

\bibitem[Martinsson(2013)]{Martinsson2013}
P.-G. Martinsson.
\newblock A direct solver for variable coefficient elliptic {PDE}s discretized
  via a composite spectral collocation method.
\newblock \emph{{J}. of {C}omp. {P}hy.}, 242:\penalty0 460--479, 2013.

\bibitem[Martinsson(2015)]{Martinsson2015}
P.-G. Martinsson.
\newblock The hierarchical {P}oincar{\'e}--{S}teklov ({HPS}) solver for
  elliptic {PDE}s: {A} tutorial.
\newblock \emph{arXiv:1506.01308}, 2015.

\bibitem[Martinsson and Rokhlin(2005)]{MartinssonRokhlin2005}
P.-G. Martinsson and V.~Rokhlin.
\newblock A fast direct solver for boundary integral equations in two
  dimensions.
\newblock \emph{{J}. of {C}omp. {P}hy.}, 205\penalty0 (1):\penalty0 1--23,
  2005.

\bibitem[Melia et~al.(2025)Melia, Fortunato, Gillman, and
  O'Neil]{MeliaFortunatoGillmanONeil2025}
D.~Melia, D.~Fortunato, A.~Gillman, and M.~O'Neil.
\newblock Hardware acceleration for hierarchical {P}oincar{\'e}--{S}teklov
  algorithms in two and three dimensions.
\newblock \emph{arXiv:2503.17535}, 2025.

\bibitem[Schmitz and Ying(2012)]{SchmitzYing2012}
P.~G. Schmitz and L.~Ying.
\newblock A fast direct solver for elliptic problems on general meshes in {2D}.
\newblock \emph{{J}. of {C}omp. {P}hy.}, 231\penalty0 (4):\penalty0 1314--1338,
  2012.

\bibitem[Xia et~al.(2009)Xia, Chandrasekaran, Gu, and
  Li]{XiaChandrasekaranGuLi2009}
J.~Xia, S.~Chandrasekaran, Y.~Gu, and X.~S. Li.
\newblock Superfast multifrontal method for large structured linear systems of
  equations.
\newblock \emph{{SIAM} {J}. on {M}atrix {A}nal. {A}ppl.}, 31\penalty0
  (3):\penalty0 1382--1411, 2009.

\end{thebibliography}

\end{document}